\documentclass[12pt]{article}

\usepackage{amssymb, latexsym} %comment
\usepackage[all]{xy}
\usepackage{times}
\newtheorem{Theorem}{Theorem}[section]
\newtheorem{Proposition}[Theorem]{Proposition}

\newtheorem{Corollary}[Theorem]{Corollary}
\newtheorem{Remark}[Theorem]{Remark}

\newtheorem{Question}[Theorem]{Question}
\newtheorem{Conjecture}[Theorem]{Conjecture}
\newcommand{\bi}{\begin{enumerate}}
\newcommand{\ei}{\end{enumerate}}
\newcommand{\be}{\begin{equation}}
\newcommand{\ee}{\end{equation}}
\newcommand{\ba}{\begin{array}}
\newcommand{\ea}{\end{array}}

\def\bih{(M,[g],J_\pm)}
\def\t{\theta}
\def\c{\mathbb{C}}
\def\cp{\mathbb{CP}}
\def\r{\mathbb{R}}
\def\j{\mathbb{J}}
\def\z{\mathbb{Z}}

\input xy
\xyoption{all}

\author{ Massimiliano Pontecorvo \\ }
\title{On a question of Vaisman concerning \\
complex surfaces}
\date{May, 2012}
\begin{document}
\maketitle

\medskip
\hfill In memory of Marco Brunella
\bigskip

\begin{abstract}
The last years have seen striking improvements on Vaisman's
question about existence of locally conformally K\"ahler (lcK) 
metrics on compact complex surfaces.
The aim of this paper is two-fold.
We review results of different authors which,
for all known examples of compact complex surfaces,
give a complete answer to Vaisman's question.
%In doing this we recall the known results as well as the open questions 
%about classification of non-K\"ahler surfaces.
\par\noindent
We also point out a relation between lcK surfaces
and generalized K\"ahler geometry 
in four-dimension and prove a new result concerning generalized K\"ahler 
structures on Hyperbolic Inoue surfaces.
We conclude with a simple observation on a question of Brunella.
\end{abstract}

%%%%                 \hfill In memory of Marco Brunella

%\bigskip
%\begin{abstract}

%The last years have seen striking improvements on Vaisman's
%question about existence of locally conformally K\"ahler metrics
%on compact complex surfaces.

%\end{abstract}

\section{Introduction}

%\bigskip
Throughout this work $M$ will denote a smooth oriented $4$-manifold -
which will be compact most of the times and always connected -
equipped with a Riemannian metric $g$, but in fact most of our
notions will only depend on the conformal class $[g]$ of all
metrics $g^\prime=e^{f}\,g$ which are obtained multiplying $\,g\,$
by a positive smooth function.

$(M,g,J)$ is said to be a Hermitian surface when $J$ is a
$g$-orthogonal endomorphism of the tangent bundle $J:TM\to TM$
with $J^2=-id.$ inducing holomorphic coordinates, so that $(M,J)$
is a smooth complex manifold and dim$_{\c}M=2$. Because $J$ is
orthogonal for any metric in $[g]$,
we will often refer to the Hermitian conformal surface $(M,[g],J)$.

\medskip

The \it Hermitian $2$-form \rm of the metric $(g,J)$ is defined to be
$\Omega(\cdot,\cdot):=g(\cdot,J\cdot)$ and in four-dimensions
defines an isomorphism between $1$-forms and $3$-forms so that
$$d\Omega =\Omega \wedge \t$$
for a unique $1$-form $\t$ called the \it Lee form \rm of the
Hermitian metric; since the metric is K\"ahler if and only if
$\t=0$ and the Lee form $\t^\prime$ of a conformally related metric
$g^\prime=e^{f}\,g \;$ is given by $\t^\prime = \t +df$, we easily have:

\bi
 \item The Lee form $\t=df$ is exact when the metric $e^{f}\,g$ is K\"ahler
 - i.e. $(g,J)$ is a globally conformally K\"ahler metric.
 \item The Lee form is closed - i.e. locally exact, $d\t=0$ -
 precisely when the Hermitian metric is conformal to K\"ahler, locally.
 \ei

\par\noindent Such metrics are then called
 \it locally conformally K\"ahler, \rm usually abbreviated by lcK.
 Of course, it is a conformally invariant notion.

\medskip
 The paper is dedicated to discuss the following question of Vaisman:

\begin{Question} {\rm \cite[p.122]{va87}}\label{vai}
Which compact complex surfaces $(M,J)$ admit lcK metrics?
\end{Question}

By a  celebrated result of Gauduchon, in the compact four-dimensional case, there
always exists a unique metric in the conformal class - called the
Gauduchon metric of $(M,[g],J)$ - whose Lee form $\t$ is coclosed.
Furthermore $\t$ is coexact if and
only if the first Betti number $b_1(M)$ of $M$ is even \cite{ga84}.

%\bigskip

As a consequence, when $b_1(M)$ is even, the Lee form of a Gauduchon
metric on $(M,J)$ can be closed only for K\"ahler metrics.
It follows that lcK compact surfaces with $b_1$ even are automatically
globally conformally K\"ahler;
this had already been  proved by Vaisman  in a different way \cite{va80}.
The conclusion is that, for $b_1(M)$ even, Vaisman question
coincides with a famous conjecture of Kodaira:

\begin{Conjecture}{\rm\cite{km}} A compact complex surface is
K\"ahlerian if and only if the first Betti number is even.
\end{Conjecture}

This  conjecture became a Theorem since the early 1980's and was first proved
using the Enriques-Kodaira classification of surfaces with the contribution of 
Miyaoka, Todorov and Siu who clarified \cite{si83} the remaining case of K3 surfaces. 
It is also interesting to notice that two independent,
direct proofs by Buchdahl and Lamari - not using surface
classification - have simultaneously appeared more recently
\cite{bu99}, \cite{la99}.

\medskip

Because any compact complex surface can be obtained by blowing up
points on a minimal model, the Enriques-Kodaira classification only
takes care of \it minimal \rm surfaces - i.e. without smooth rational
curves of self-intersection $(-1)$ - and fits very well with the following
result of Tricerri which generalizes the analogous well known result
of Kodaira in the K\"ahler case; see also the work of Vuletuscu
for a more recent reference:

\begin{Proposition}{\rm \cite{tr82}, \cite{vu09}}  A complex manifold $M$ is lcK
if and only if the blow up of $M$ at a point is lcK
\end{Proposition}

As noticed in \cite[Remark 4.3]{tr82} this reduces Vaisman
question to minimal surfaces and it is quite natural to use the
Enriques-Kodaira classification to attack the problem.

\bigskip
For the reader's convenience we will display  - at the end of this section - 
a (slightly refined) table \cite[p.244]{bpv} of this well known classification which divides
all minimal compact complex surfaces $S:=(M.J)$ in eleven classes
belonging to four groups according to the possible values of the \it
Kodaira dimension \rm $Kod(S)=-\infty,0,1,2,$ which appears in the
second column of the table; in the other columns are indicated:
the \it algebraic dimension \rm $a(S)$, which is the degree of the
extension (over $\c$) of the field of meromorphic functions on $M$,
the Euler characteristic $\chi(S)=2-2b_1(S)+b_2(S)$ and the first
Betti number $b_1(S)$.

\medskip
We believe that Vaisman question is a good one
because of its intricate relations with the study of complex geometry of non-K\"ahler  surfaces
developed by the Japanese school of Kodaira, Nakamura, Inoue, Kato, Enoki
and more recently by the Marseille-based group of Dloussky, Teleman, Oeljeklaus, Toma.

\bigskip
The remaining of the paper is divided in three sections;
the first is dedicated to the work of Belgun \cite{be00} who
%- in his thesis under Gauduchon -
gave new impetus to the conjecture by completing the work initiated by Vaisman, Tricerri and Gauduchon-Ornea \cite{go} on complex surfaces of vanishing Euler characteristic.
The case of class-VII surfaces appears to be the most difficult and this
parallels the fact that the classification of complex structures on these surfaces ($\chi =0$)
was solved by Bogomolov only in 1976 \cite{bg76}.
The proof was later clarified by Li-Yau-Zheng \cite{lyz} and Teleman \cite{te94},  independently.

One major discovery of Belgun is that - unlike in the K\"ahler case - lcK is not an open
property under holomorphic deformations.
More precisely, he showed that some of the Inoue-Bombieri surfaces \cite{in74}, \cite{bm73}
are not lcK.
As a matter of fact this Hermitian property is somehow reflected in the complex geometry:
Inoue-Bombieri surfaces are the only known compact complex surfaces with no curves.

\medskip
After the work of Belgun, Vaisman question had a complete (although not always positive) 
answer for all surfaces except for those with minimal model satisfying $b_1=1$ and $b_2 >0$;
they are usually said to belong to class VII$_0^+$.
In section 3, we present our contribution about Inoue surfaces in this class.

It is important at this point to spend a few words on the classification problem 
of class VII$_0^+$ surfaces,
the best reference for known results and open questions
are the expositions of  Nakamura \cite{na89}, \cite{na08}.
It is not even known - for example -  whether every $S\in $VII$_0^+$ has a (holomorphic) curve;
a remarkable, recent result of A.Teleman \cite{te10} implies that this is the case
when $b_2(S)\leq 2$.

As a matter of fact, the only known examples in this class are so called Kato surfaces
\cite{ka77}:
they can be characterized by saying that they have $b_2(S)$ rational curves \cite{dot}.
Among them, the Inoue surfaces \cite{in77} form a subclass in which
every rational curve belongs to a cycle  $C\subset S$.
%\smallskip
These Inoue surfaces where the first examples of lcK surfaces in class VII$^+_0$.
Their lcK metrics came from a link - pointed out by Boyer,
see Prop. \ref{lck} - with self-duality in dimension four.

Surprisingly enough, it turns out that many Inoue surfaces in class VII$_0^+$
carry bi-Hermitian metrics.
This provides a new relation between lcK and generalized K\"ahler geometry  
\cite{gu07}, \cite{hi06} given by the fact that anti-self-dual bi-Hermitian metrics always satisfy 
the generalized K\"ahler property that the sum of the two Lee forms vanishes identically
\cite[3.11]{po97}.

We present a new result in this context; Theorem \ref{hyp} proves existence of generalized K\"ahler
metrics on all admissible blow-ups of Hyperbolic Inoue surfaces.
The first lcK examples in class VII$_0^+$ are by LeBrun \cite{le91} on Parabolic Inoue surfaces;  
later, a twistor construction by Fujiki-Pontecorvo \cite{fp10} produced these very special 
lcK metrics on all (minimal) Hyperbolic as well as on all Half-Inoue surfaces.
%The fact that these examples admit deformations which are again lcK seemed to indicate that indeed
%many Kato surfaces are lcK.

\medskip
Shortly after came the remarkable contribution of Marco Brunella \cite{br10}, \cite{br11},
which is the subject of the last section.
He proved that - contrary to the case of Inoue-Bombieri surfaces -
the lcK property is open under holomorphic deformations of Kato surfaces.
In fact, he showed furthermore that all Kato surfaces are lcK

Once again, Brunella results reflect the classification problem of 
non-K\"ahler surfaces: indeed, Kato surfaces are the only known surfaces in class VII$_0^+$.

Vaisman question would now have a complete answer if no other class-VII$_0^+$
surfaces would exist. This is in fact a strong conjecture of Nakamura  \cite[5.5]{na89}.

Much more recently  A.Teleman, using Gauge theory, has proved this conjecture for 
$b_2=1$ in \cite{te05};
he also proved in \cite{te10} that for $b_2=2$ every calss-VII$_0^+$ surface 
contains a cycle of rational curves.

We conclude the paper with a simple observation on a question of Brunella \cite{br11}.

\bigskip

\noindent \bf Acknowledgement. \rm
A preliminary and incomplete version of this work
appeared in the informal volume on the occasion
of Fujiki's 60th birthday \cite{po08}.
We take this opportunity to thank the organizers of the conference.

\newpage

\begin{large}

\begin{center} Table of Enriques-Kodaira classification
\end{center}

\end{large}
%\begin{large}

%\medskip
\begin{center} %\label{table}
\begin{tabular} {|p{7cm}|c|c|c|c|}
\hline %\vskip{.1cm}
&&&&\\
Class of $S$                  &$Kod(S)$ &$a(S)$&$\chi(S)$&$b_1(S)$ \\
&&&&\\  \hline  &&&&\\
1) rational surfaces          &               & 2    &  3,4    &  0    \\
2) ruled surfaces of genus $g\geq 1$&$-\infty$& 2    &$4(1-g)$ &  2g \\
3) class $VII_0$ surfaces     &               & 0,1  &  0      &  1    \\
4) class $VII_0^+$ surfaces   &               & 0    &  $>0$   &  1    \\
&&&&\\    \hline    &&&&\\
5) tori                       &         & 0,1,2& 0       & 4     \\
6) K3-surfaces                &         & 0,1,2& 24      & 0     \\
7) hyperelliptic surfaces     &   0     & 2    & 0       & 2     \\
8) Enriques surfaces          &         & 2    & 12      & 0     \\
9) Kodaira surfaces           &         & 1    & 0       & 1,3   \\
&&&& \\ \hline &&&& \\
10) properly elliptic surfaces & 1       & 1,2    &$\geq 0$ & even  \\
                               &         & 1      & 0       & odd   \\
&&&& \\ \hline &&&& \\
11) surfaces of general type  & 2       & 2    & $>0$    & even  \\
&&&& \\ \hline
\end{tabular}
\end{center}

%%%%%%%%%%%%%%%%

%\bigskip

\section{%The work of Belgun: non-K\"ahler surfaces with $\chi=0$ or
               Non-K\"ahler surfaces with $\chi=0$}

After reducing the problem to minimal surfaces with $b_1(S)$
odd, we see from Kodaira classification that the Euler characteristic
satisfies $\chi(S)\geq 0$. Extending previous results of Vaisman
\cite{va87}, Tricerri \cite{tr82} and Gauduchon-Ornea \cite{go},
Belgun gave a complete answer to Vaisman question in the limiting
case when $\chi(S)=0$, see \cite{be00}.
 In fact, because $\chi(S)=0$ is clearly a necessary condition for a
lcK metric to have parallel Lee form, Belgun even classified
compact complex surfaces with these special lcK metrics which
are called \underline{Vaisman metrics}.
The main tool was the use of explicit models for the universal cover
of such geometric complex surfaces, which we now describe in order of
decreasing Kodaira dimension. Recall that in what follows we assume
$b_1(S)$ odd and $\chi(S)=0$ which implies $S$ minimal, by Kodaira
classification and the fact that $\chi$ increases under blow up.

\medskip\noindent
\underline{$Kod(S)=1$}\par\noindent
\smallskip
These are so called properly elliptic surfaces which means that they
are elliptic of algebraic dimension 1. In general, $Kod(S)\leq a(S)$
and in our case equality holds because every Moischezon surface is
algebraic \cite{bpv}. Then $S$ has a holomorphic map to a curve with
connected fibers and the generic fiber is elliptic, otherwise it
would have non-zero intersection with the canonical class and $S$
would be algebraic.  %FM 132
 Elliptic surfaces were studied in depth by Kodaira \cite{ko64} and
when $b_1$ is odd, the only singular fibers turn out to be multiple
fibers. Maheara \cite{ma77} showed that the universal cover is
$\c\times H$ the product of a complex line with the upper half
plane. The explicit description of the covering map allows Belgun
to prove that all minimal surfaces with odd $b_1$ and $Kod(S)=1$
admit lcK metrics and in fact Vaisman metrics.

\medskip\noindent
\underline{$Kod(S)=0$}
\smallskip \par\noindent
These are called Kodaira surfaces and are divided into primary or
secondary according to whether $b_1$ is $1$ or $3$. Primary
Kodaira surfaces $S=(M,J)$ are elliptic fiber bundles over an
elliptic curve.
It was shown by Thurston that they also admit a symplectic
structure $I$ and $(M,I)$ provided the first example of 
a compact symplectic manifold without K\"ahler stuctures \cite{th76}.
%\cite{ab84}
Finally, Salamon observed in \cite{sa94} that $I$ and $J$ anti commute
and therefore give rise to an almost hyperHermitian structure spanned
by $\{I,J,IJ\}$ with one complex structure $J$ and two symplectic structures $I$ and $IJ$.

A secondary Kodaira surface is finitely covered by a primary one and the universal
covering is a nilpotent group; Belgun showed that all of these surfaces
always admit Vaisman metrics, in particular lcK metrics.

\medskip\noindent
\underline{$Kod(S)=-\infty$}
\smallskip \par\noindent
This is the class of VII$_0$-surfaces in the above table, they are
characterized by having negative Kodaira dimension and vanishing
second Betti number; the already cited Bogomolov theorem \cite{bg76} 
states they are either Hopf surfaces - which always have at least one curve - 
or the surfaces independently discovered at the same time by Inoue
\cite{in74} and Bombieri \cite{bm73} which have no curves at all;
we call them Inoue-Bombieri surfaces in order to distinguish from 
a different type of Inoue surfaces \cite{in77} which belong to class VII$_0^+$
and will be the subject of next section.

\underline{Hopf surfaces} have universal covering $\c^2\setminus 0$
on which we take global coordinates $(z,w)$; the primary Hopf
surfaces are diffeomorphic to $S^1\times S^3$ with fundamental group
generated by the contraction
$$(z,w)\mapsto (az,bw+\lambda z^n)$$
where $a,b,\lambda\in\c$ and $n\in {\mathbb N}$ satisfy
$0<|a|<|b|<1$ and $\lambda(a-b^n)=0$; we say that a Hopf surface
is \it diagonal \rm when $\lambda =0$ and in this case there are
always (at least) two (elliptic) curves on  $S$, namely the images
of the two axes $\{z=0\}$ and $\{w=0\}$. When $\lambda\neq 0$ the
only curve on $S$ is the image of $z=0$. Every Hopf surface is
finitely covered by a primary one.

Extending results of \cite{go} Belgun proved that every Hopf surface
is lcK; furthermore, the Lee form can be parallel exactly for the diagonal
ones.

\smallskip
\underline{Inoue-Bombieri surfaces} come in three different
families which for simplicity we denote by $S_M$, $S^-_N$ and
$S^+_{N,u}$ with $u\in\c$; the universal covering is always a
solvable group biholomorphic to $\c\times H$.

Tricerri \cite{tr82} constructed explicit lcK metrics on every surface of
type $S_M$ and $S^-_N$, as well on the surfaces
$S^+_{N,u}$ \underline{with $u\in\r$}.
These where also the first lcK metrics with non-parallel Lee form.
Belgun showed that these examples are in fact sharp.

\smallskip
\par\noindent Finally, we state the complete result

\begin{Theorem}{\rm\cite{be00}} \label{florin} 
 (1) All diagonal Hopf surfaces and all surfaces with
 $Kod\geq 0$ admit lcK metrics with parallel Lee form.
 \par\noindent (2) All other Hopf surfaces and the Inoue-Bombieri surfaces $S_M$, $S_N^-$ and
 $S^+_{N,u}$ with $u\in\r$, admit lcK metrics whose Lee form cannot be parallel.
 \par\noindent (3) When $u\in\c\setminus\r$, the Inoue-Bombieri surfaces $S^+_{N,u}$
 do not admit lcK metrics at all.

 \smallskip
\noindent Furthermore, by Kodaira classification, every compact complex
surface with $b_1$ odd and $\chi=0$ belongs to exactly one of the
above three classes.
\end{Theorem}

In particular, contrary to the K\"ahler case \cite{ks60}, the family of
Inoue-Bombieri surfaces $S^+_{N,u}$ shows that

\begin{Corollary}{\rm\cite{be00}}
The class of lcK manifolds is not open under holomorphic deformations.
\end{Corollary}

%%%%%%%%%%%%%%%%%%%%%%%%%

\section{%Anti-self-dual bi-Hermitian metrics or
               Inoue surfaces in class VII$_0^+$}

After the work of Belgun, Vaisman question had a complete (although
not always positive) answer for all surfaces except for the class of
minimal surfaces with $b_1$ odd and $\chi\neq 0$.
We see from Kodaira classification that this in fact implies:
$b_1=1$,  $\chi = b_2 = b_2^+ > 0$ and $Kod = -\infty$.
This class of surfaces is usually denoted by VII$_0^+$.

The subclass with $b_2(S)$ rational curves each of which belongs to a cycle $C$
are called Enoki \cite{en81} and Inoue surfaces \cite{in77}; notice that
there can be at most two cycles of rational curves on a given surface \cite[6.3]{na89}.

%Furthermore, there are four types of Inoue surfaces in class VII$_0^+$: they are called
Enoki surfaces  are characterized by the topological property 
that the unique cycle satisfies $C^2=0$;  some of the Enoki surfaces 
also admit an elliptic curve in which case are called Parabolic Inoue surfaces.
There are two more types: Hyperbolic Inoue surfaces, have two cycles of rational curves
and Half-Inoue surfaces which have a unique cycle with $C^2<0$;
the name comes from the fact that they are double-covered by Hyperbolic Inoue surfaces.

\bigskip

The first examples of lcK complex surfaces in class-VII$_0^+$ where found
by LeBrun \cite[p.391]{le91} on some Parabolic Inoue surfaces;
this is because of the following link %which was
pointed out by Boyer, see also \cite{po91}.

\begin{Proposition}{\rm\cite{bo86}}\label{lck} Hermitian
anti-self-dual metrics on compact four-manifolds are necessarily lcK
\end{Proposition}

It was later noticed \cite{po97} that LeBrun's examples are in fact {\it bi-Hermitian}.
This means that we have a four-manifold $M$ with two integrable complex structures $J_\pm$ 
which are orthogonal w.r.t. a Riemannian metric and induce the {\it same orientation}.
Of course we also require that there is a point $p\in M$ where $J_-(p)\neq\pm J_+(p)$;
this determines the conformal structure uniquely and is equivalent to ask that the following
algebraic condition holds \cite{tr75}
$$J_+J_- + J_-J_+ = -2a\cdot id$$
for some smooth $a:M\to [-1,1]$ usually called the angle function \cite{po97}.
Notice that this is a conformally-invariant notion, for which we will use the notation $\bih$.
In order to understand the structure of such surfaces we now give a twistor presentation.

It is sometimes useful
to think of an almost Hermitian structure $J$ as a smooth section
$J:M\to Z$ of the twistor space. $Z$ is the fiber bundle of all linear
complex structures at $T_p M$ compatible with the metric and orientation;
in four-dimension the fiber at $p\in M$ is the homogeneous space $SO(4)/U(2)\cong \cp_1$.
It is known that $Z$ is an almost complex $6$-manifold - let $\j$ denote
the almost complex structure - which only depends on the fixed
conformal structure $[g]$ and orientation of $M$. The integrability
of $J$ is equivalent to the fact that $J:M\to Z$ is an almost
$\j$-holomorphic map and is also equivalent to ask that its image
$J(M)=:S$ is an almost complex submanifold of $(Z,\j)$, see for example \cite{amp};
$S$ is therefore (tautologically) biholomorphic to the original complex surface:
$S\cong (M,J)$.

The twistor space $Z$ is equipped with a natural \it real
structure, \rm namely the involution $\tau:Z\to Z$ that sends
$J\mapsto -J$ and objects that are $\tau$-invariant are then called
`real'. For example each fiber $SO(4)/U(2)\cong \cp_1$ is called a
`real twistor line' because it is a $\j$-complex curve
biholomorphic to $\cp_1$ and defines a non-trivial cohomology
class $h\in H^2(Z,\z)$. It is then known that the first Chern
class $c_1(Z,\j)=-4h$; when $M$ is compact we also have another
`real' 4-cohomology class if we let $X$ denote the (obviously
disjoint) union $S \amalg \tau(S)$.
 Then its Poincar\'e dual $X^*$ equals $2h$ because $X$ intersects each
real twistor line in exactly two points with the complex
orientation and we can now compute the first Chern class of
$(M,J)\cong S$ from adjunction formula:
 \be \label{adj} c_1(M,J)=[c_1(Z)+ X^*]_{|_{S}}=-2h_{|_{S}} \ee

In the bi-Hermitian case there are two sections $J_\pm:M\to Z$;
their images define two complex hypersurfaces
$S_\pm:=J_\pm(M)$ and we will also consider the real hypersurfaces
$X_\pm:=S_\pm\amalg \tau(S_\pm)$. Notice again that there is a
tautological biholomorphism $(M,J_\pm)\cong S_\pm$.

Next, we concentrate on the following intersection

$$C_+:=S_+\cap X_- \subset Z$$

We can think of $C_+$ as a divisor and therefore a holomorphic line bundle on $S_+\cong (M,J_+)$:
it is a complex 1-dimensional subset because it is the intersection of (three) smooth complex
hypersurfaces (with no common component) of the almost complex manifold $(Z,\j)$.
The multiplicities of the irreducible components of $C_+$ are the (necessarily positive)
intersection multiplicities and are in any case uniquely determined by the Chern
class which we compute below.

%Being the intersection of (three) smooth complex hyper surfaces of the
%almost complex manifold $(Z,\j)$, $C_+$ is an effective (or empty) divisor
%in the complex surface $S_+\cong (M,J_+)$.

In particular, its image $C:=t(C_+) \subset M$ - which is easily seen to be the set of points where
$J_\pm$ are linearly dependent - is a closed nowhere dense subset of $M$  \cite[1.3]{po97}
as well as the support of the divisor $C_+$.

$$C=t(S_+\cap S_-)\amalg t(S_+\cap \tau(S_-))=\{p\in M\, |\, J_+(p)=\pm J_-(p)\}$$

The case $C=\emptyset$ is
usually referred to as \it strongly bi-Hermitian \rm because
$J_\pm$ are linearly independent everywhere on $M$ in this
situation; see \cite{ad} for a complete classification.

We now proceed to compute the first Chern class of the divisor $C_+$ in $S_+$.

\begin{Proposition} The first Chern class of $C_+$ is the anti-canonical class of
the surface $S_+$:
 $$c_1(C_+)=-c_1(S_+)$$
 i.e. on the compact complex surface $(S_+)$ there is a holomorphic
line bundle $L_+$ such that

 (i) \quad $c_1(L_+)=0$, and

 (ii) \quad $C_+=L_+ - K$ where $K$ is the canonical line
 bundle of the surface $S_+$.
\end{Proposition}

\it Proof. \rm In the twistor space $Z$ we have that
 $C_+\cong X_-\cap S_+$ therefore $c_1(C_+)=2h_{|{S_+}}$
and the conclusion follows from (\ref{adj}). \hfill $\Box$

\begin{Remark} In other words $C_+$ is a numerically anti-canonical
divisor in the complex surface $S_+$. By a powerful result of Dloussky
{\rm \cite{dl06}} this forces a class-VII$_0^+$ surface to actually be a Kato surface.
We also recall that a numerically anti-canonical divisor is
automatically effective on a Kato surface {\rm \cite[2.2]{dl06}}.

Of course we also have a numerically anti-canonical divisor $C_-$
in the other complex surface $S_-\cong (M,J_-)$ with the same support $t(C_-)=t(C_+)=C$.
This implies that $(M,J_\pm)$ are both Kato surfaces;
although they don't need to be biholomorphic in general {\rm \cite[7.5]{fp10}}
the fact that both have a numerically anti-canonical divisor with the same irreducible
components and multiplicities is enough to conclude that their complex
structures share a lot of common properties.
\end{Remark}

\bigskip
Perhaps surprisingly, the above presentation is not only useful to describe
the complex structure of a possible bi-Hermitian surface but it proved to be
effective in producing new examples.
Starting from a Joyce twistor space \cite{jo95} \cite{fu00} over the connected sum
$m\mathbb{CP}_2$ we produce in \cite{fp10} the twistor space of bi-Hermitian anti-self-dual
metrics on surfaces of class VII$_0^+$ with $b_2=m$.

It is known that this kind of metrics can only leave on Parabolic or Hyperbolic
Inoue surfaces and our result covers both cases;
%by \ref{lck} they are automatically lcK. In particular, we proved the following
furthermore, recall that the following strong form of \ref{lck} applies:

\begin{Proposition}{\rm \cite[3.1, 3.11]{po97}} \label{bih}
Let  $(M,[g],J_\pm)$ be a compact bi-Hermitian surface.
Then, the metric $[g]$ is anti-self-dual if and only if $([g],J_-)$ is lcK
if and only if $([g],J_+)$ is lcK; in this situation the $J_+$-Gauduchon metric
is automatically $J_-$-Gauduchon and there are only two cases:

\noindent (i) The two Lee forms coincide: $\t_+=\t_-$ and $(g,J_\pm)$ span a 
  hyper-Hermitianstructure, or else

\noindent (ii) The sum of the two Lee forms vanishes identically: $\t_+ + \t_- =0$ 
\end{Proposition}

Our twistor construction produced in particular the following result:

\begin{Theorem}{\rm \cite[7.5, 9.3]{fp10}}
Every Hyperbolic Inoue surface admits anti-self-dual bi-Hermitian structures.
Every Half Inoue surface admits anti-self-dual Hermitian structures.
By \ref{lck} and \ref{bih} all such surfaces are lcK.
\end{Theorem}

It is also interesting to point out that bi-Hemitian metrics provide a link between lcK
and generalized K\"ahler geometry in four-dimension.
This is because a generalized K\"ahler structure induces a bi-Hermitian structure 
satisfying the second equation in \ref{bih}, see \cite[Prop. 3]{ad}.
%This condition - already introduced in \cite{agg} -
%is satisfied by anti-self-dual bi-Hermitian metrics \cite{po97}[3.11],
%provided they are not hyper-Hermitian (in which case the two Lee forms coincide).
We can now prove the following optimal and new result about
%bi-Hermitian metrics on
non-minimal Hyperbolic Inoue surfaces.

\begin{Theorem}  \label{hyp}
The blow up of any Hyperbolic Inoue surface at any set of points
on its anti-canonical divisor
admits bi-Hermitian metrics which are generalized K\"ahler.
\end{Theorem}
{\em Proof.} The already mentioned result \cite[7.5]{fp10} 
produces in fact anti-self-dual bi-Hermitian metrics
on any \it properly blown-up \rm Hyperbolic Inoue surface - which means that
the blown-up points are singular points of the anti-canonical divisor.
Since the only compact hyper-Hermitian surfaces are Hopf surfaces \cite{bo88},
these examples must satisfy equation (ii) in \ref{bih} and are therefore generalized K\"ahler.
We can now apply a recent ``complementary" result of Cavalvanti-Gualtieri \cite{cg11} 
which states that one can blow up smooth points of the anti-canonical divisor and
still have a generalized K\"ahler structure.
\hfill$\Box$

\begin{Remark} The condition that the blown up points belong to the
anti-canonical divisor is necessary {\rm \cite[3.14]{po97}}.
\end{Remark}

%We conclude this section with the observation that our bi-Hermitian examples coming from twistor
%theory admit large families of deformations which are automatically lcK.
%This can be seen as an indication that ``many" Kato surfaces are lcK.

%%%%%%%%%%%%%%%%

%\bigskip
\section{ %The work of Brunella or
   Kato surfaces}

After the work of LeBrun and Fujiki-Pontecorvo
every known example of lcK surfaces in class VII$_0^+$
where Inoue surfaces (a particular case of Kato surfaces)
and the metrics were actually anti-self-dual
and, up to double coverings, bi-Hermitian.

\smallskip
This situation was clarified by Marco Brunella, who left us much too early.
In a sequence of two papers he proved that Kato surfaces behave extremely
well with respect to the lcK property.
We like to think that the results of the previous section
somehow inspired him to work on the subject.
This section starts with some preliminaries about Kato surfaces \cite{ka77}.

\smallskip
%All other Kato surfaces have branches: there is a unique cycle $C$ but not all rational
%curves belong to $C$.

Surfaces in class VII$_0^+$ all have vanishing algebraic dimension
and therefore at most finitely many (holomorphic) curves.
All known examples in this class are so called \it Kato
surfaces: \rm they contain a global spherical shell which means that
they are obtained by holomorphic surgery from a neighborhood of the
origin in $\c^2$ blown up at ``infinitely-near" points, in
particular they are diffeomorphic to $(S^1\times S^3)\# m\overline{\cp}_2$ where $m=b_2(S)$.

VII$_0^+$-surfaces have been intensively studied in the 1980's  by the
Japanese school of Nakamura, Kato, Inoue, Enoki; they proved that
every Kato surface has small deformations which are blown-up Hopf
surfaces. Nakamura \cite{na90} also showed that every Kato surface
has a singular degeneration which is a rational surface with a nodal
curve; the case $b_2=0$ actually goes back to Kodaira \cite{ko64}.
This characterization very much inspired our work \cite{fp10} because
we use a twistorial version of Nakamura degeneration.

There are several open problems concerning the complex structure
of a surface $S\in $VII$_0^+$ the basic questions being: does $S$ admit curves?
Or more strongly: does it admit a cycle of rational curves? Or even more: is $S$
necessarily a Kato surface?

Recent progress has been made by the Marseille-based group of
Dloussky, Oljeklaus, Toma, Teleman;
by \cite{dot} a Kato surface is characterized in class VII$_0^+$ by having $b_2$ rational curves
or equivalently,  a numerically anti-canonical divisor \cite{dl06}.
More recently A.Teleman has shown the following result by studying the moduli space of $PU(2)$ instantons

\begin{Theorem}{\rm \cite{te10}}
Every class VII$_0^+$ with $b_2\leq 2$ has a cycle of rational curves.
In particular, is diffeomeorphic to $(S^1\times S^3) \# b_2 \overline{\cp}_2$ and
must be a Kato surface when $b_2=1$.
\end{Theorem}

The first results about Hermitian geometry of these surfaces where the ones of the
previous sections.
Thanks to the work of Brunella this situation changed dramatically soon after.
His first result is that Kato surfaces are well-behaved with respect to deformations
of lcK metrics.

\begin{Theorem} {\rm \cite{br10}}
For Kato surfaces, the property of being lcK is open under holomorphic deformations.
\end{Theorem}

As a consequence, Bombieri-Inoue surfaces of type $S^+_{N,u}$ are the only known compact complex surfaces which are not lcK-stable by small deformations.

\medskip

The above Theorem, combined with the results of the previous section, already gives
many new lcK surfaces. However, it was greatly generalized in a following work

\begin{Theorem} {\rm\cite{br11}} \label{marco} 
Every Kato surface admits lcK metrics.
\end{Theorem}

This is a remarkable result, the proof is short and smart;
it involves a careful construction of K\"ahler potentials on the modified unit ball
in such a way that the holomorphic surgery becomes actually a conformal map.

\bigskip
To summarize, we state the results described in this paper in the following

\begin{Remark}
Theorems \ref{florin} and \ref{marco} answer Vaisman
question \ref{vai} for all known compact complex surfaces.

Furthermore, Nakamura strong Conjecture {\rm \cite[5.5]{na89}}
that every class VII$_0^+$ is actually a Kato surface,  would then imply
that every compact complex surface is lcK except for Bombieri-Inoue surfaces
$S_{N,u}^+$ with $u\in\c\setminus\r$, which are not lcK
\end{Remark}

We finish this work with a simple observation on \cite[question after Corollary 2]{br11}.

\begin{Remark}
The universal covering is K\"ahlerian for every known compact complex surface.

In fact, a lcK metric is the same as a K\"ahler metric on the universal cover which
transforms by homotheties under the action of the group of deck transformations.
Furthermore,  Inoue-Bombieri surfaces are covered by 
the product of the upper half-plane with the complex line {\rm \cite{in74}}.
Therefore, the universal cover happens to be K\"ahlerian even for the
non lcK surfaces of Belgun.
\end{Remark}

%In fact, this certainly holds true for the compact lcK surfaces, but it
%happens to apply to the universal covering of the non-lcK surfaces $S_{N,u}^+$
%as well; just because any Bombieri-Inoue surface is covered by
%the product of the upper half-plane with the complex line \cite{in}.

%%%%%%%%%%%%%%%%%

\newcommand{\bysame}{\leavevmode\hbox to3em{\hrulefill}\,}

\bigskip \noindent
M.Pontecorvo
\par\noindent Dipartimento di Matematica
\par\noindent Universit\`a Roma Tre
\par\noindent L.go S.L.Murialdo 1, 00146 Roma  
\par\noindent max@mat.uniroma3.it
\end{document}